\documentstyle[12pt]{article}

\def\nl#1#2{\mathop{#2}\limits_{#1}}

\newcommand{\pl}{\partial}

\newcommand{\ol}{\overline}
\newcommand{\bqq}{\begin{equation} \label}
\newcommand{\eeq}{\end{equation}}

\begin{document}

\begin{center}
{\Large{\bf{Metrics of 6-dimensional $h$-spaces}}}\\
{\Large{\bf{of the $[3(21)]$, $[(32)1]$, $[(321)]$ types
}}}
\end{center}

\bigskip

\centerline{{\large
Zolfira Zakirova}\footnote{e-mail: Zolfira.Zakirova@soros.ksu.ru, Zolfira.Zakirova@ksu.ru}}

\bigskip

\bigskip

\abstract{\small  In this note we find the metrics of 6-dimensional
$h$-spaces of the $[3(21)]$, $[(32)1]$, $[(321)]$ types and then determine the
quadratic first integrals of the geodesic equations of these $h$-spaces.
}

\begin{center}
\rule{5cm}{1pt}
\end{center}

\bigskip
\setcounter{footnote}{0}
\section{Introduction}

The aim of this paper is to investigate the $6$-dimensional pseudo-Riemannian space $V^6(g_{ij})$
with signature
$[++----]$. In particular, we find the metrics of 6-dimensional
$h$-spaces of the $[3(21)]$, $[(32)1]$, $[(321)]$ types and then determine the quadratic first integrals
of the geodesic equations of these $h$-spaces.

The general method of determining pseudo-Riemannian manifolds that
admit some nonhomothetic projective group $G_r$ has been developed
by A.V.Aminova \cite{am2}.

It is known that the differentiable symmetric bilinear form
$h$ on the pseudo-Riemannian manifold $M^n$ satisfies the
Eisenhart equation provided the condition
\bqq{1}
\nabla h(X, Z, W)=2g(X, Z)W\varphi+2g(X, W)Z\varphi+2g(Z, W)X\varphi
\eeq
is satisfied for a $0$-form $\varphi$ in $M^n$ and arbitrary
vector fields $X, Z, W \in TM^n$.

After changing variables
\bqq{2}
h=a+2\varphi g,
\eeq
$a$ being a symmetric bilinear form of the same
characteristics as $h$, equation (\ref{1}) is rewritten
in the form \cite{am2}
\bqq{3}
\nabla a(X, Z, W)=g(X, W)Z \varphi+g(W, Z)X \varphi,
\eeq
$X, W, Z$ being arbitrary vector fields, defined within the
domain $V$.

\section{Eisenhart equation in the skew-normal frames\\
for $h$--spaces of the $[3(21)]$, $[(32)1]$, $[(321)]$ types}

\bigskip

If $\lbrace X_l \rbrace$ are skew-normal frames, then there is the following equation \cite{am2}
\bqq{4}
X_r \overline a_{pq}+\sum_{h=1}^n e_h( \overline a_{hq} \gamma _{\tilde hpr}+
\overline a_{ph} \gamma_{\tilde  hqr})=\overline g_{pr}X_q \varphi+\overline g_{qr}X_p \varphi
\eeq
$$
(p, q, r, \tilde h=1,\ldots,n),
$$
equivalent to equation (\ref{1}). Here
$$
X_r \varphi \equiv {\nl{r}\xi}^i \frac{\partial \varphi}
{\partial x^i},
\quad
\gamma_{pqr}=-\gamma_{qpr}=
{\nl{p}{\xi}\hspace{-2.5mm}\phantom{a}_{i,j}}{\nl{q}\xi}^i{\nl{r}\xi}^j,
$$
${\nl{i}\xi}^j$ are the components of skew-normal frames, $\overline a_{pq}$ and
$\overline g_{pr}$ are the canonical forms of the tensors $a_{ij}$, $g_{ij}$, respectively.

In the considered $h$-spaces the canonical forms are of the following form \cite{pet}
\bqq{5}
\ol{g}_{ij} dx^i dx^j=e_3 (2 dx^1 dx^3+{dx^2}^2) +2 e_5 dx^4 dx^5+e_6 {dx^6}^2,
\eeq
$$
\ol{a}_{ij} dx^i dx^j=e_3 \lambda_3(2 dx^1 dx^3+{dx^2}^2)+2 e_3 dx^2 dx^3+2e_5 \lambda_5 dx^4 dx^5+
$$
$$
e_5{dx^5}^2+ e_6  \lambda_{6} {dx^6}^2,
$$
$$
(e_1=e_2=e_3, e_4=e_5),
\quad
e_i =\pm 1,
$$
where $\lambda_1=\lambda_2=\lambda_3$, $\lambda_4=\lambda_5$, $\lambda_6$ are the roots of the
characteristic equation ${\det}(h_{ij}-{\lambda}g_{ij})=0$.
In the case of $h$-space of the type $[3(21)]$ $\lambda_4=\lambda_5=\lambda_6$,
for $h$-space of the type $[(32)1)]$
$\lambda_1=\lambda_2=\lambda_3=\lambda_4=\lambda_5$, in the case of $h$-space of the type $[(321)]$
all $\lambda_i$ $(i=1,\ldots,6)$ coincide.

Substituting in (\ref{4}) the canonical forms $\ol{g}_{pq}$ and  $\ol{a}_{pq}$
from (\ref{5}) and taking into account that,
for the types $[3(21)]$, $[(32)1]$, $[(321)]$ $h$-spaces,
$\tilde1=3, \; \tilde2=2, \; \tilde3=1, \; \tilde4=5, \; \tilde5=4, \; \tilde6=6$,
one gets the system of equations, which after transformations can be
written in the form

\bigskip

\noindent
characteristic $\chi_{10}=[3(21)]$:
\bqq{6}
X_r \lambda_3=0 \; (r \ne 3),
\quad
X_r \lambda_6=0,
\quad
X_3( \lambda_3-2/3 \varphi)=0,
\eeq
$$
\gamma_{213}=\frac{1}{2}e_3 X_3 \varphi,
\quad
\gamma_{312}=\gamma_{321}=-e_3 X_3 \varphi,
$$
$$
\gamma_{345}=\gamma_{354}=\frac{e_5  X_3 \varphi}{\lambda_3-\lambda_5},
\quad
\gamma_{366}=\frac{e_6  X_3 \varphi}{\lambda_3-\lambda_6},
$$
$\gamma_{56r}$ are arbitrary, all other $\gamma_{pqr}$ are equal to zero.

\bigskip

\noindent
characteristic $\chi_{11}=[(32)1]$:
\bqq{7}
X_r \lambda_5=0,
\quad
X_r \lambda_6=0 \; (r \ne 6),
\quad
X_6 (\lambda_6- 2\varphi)=0,
\eeq
$$
\gamma_{163}=\gamma_{262}= \gamma_{361}= \frac{e_3 X_6 {\varphi}}{\lambda_5-\lambda_6},
\quad
\gamma_{263}=\gamma_{362}= -\frac{e_3 X_6 \varphi}{(\lambda_5-\lambda_6)^2},
$$
$$
\gamma_{363}=\frac{e_3 X_6 \varphi}{(\lambda_5-\lambda_6)^3},
\quad
\gamma_{564}= \gamma_{465}= \frac{e_5 X_6 {\varphi}}{\lambda_5-\lambda_6},
\quad
\gamma_{565}=- \frac{e_5 X_6 \varphi}{(\lambda_5-\lambda_6)^2},
$$
$$
\gamma_{25r}=\gamma_{34r},
$$
$\gamma_{35r}$ are arbitrary, all other $\gamma_{pqr}$ are equal to zero.

\bigskip

\noindent
characteristic $\chi_{12}=[(321)]$:
\bqq{8}
X_r \lambda_6=0,
\quad
\gamma_{25r}=\gamma_{34r},
\eeq
$\gamma_{ksr}$ $(k, s=3, 5, 6, \; k\ne s)$ are arbitrary, all other $\gamma_{pqr}$  are equal to zero.

\section{$h$--space of the type $[3(21)]$}

The linear system of PDE's
$$
X_q \theta={\nl{q}\xi}^i\pl_i \theta=0,
\quad
(q=1,\ldots,m, i=1,\ldots,6),
$$
where ${\nl{q}\xi}^i$ are components of the frames above,
is known to be integrable, i.e. to admit $6-m$ independent solutions,
if all the commutators of the system (\cite{ezen1}, see also \cite{am2})
\bqq{9}
(X_q, X_r) \theta=X_q X_r \theta-X_r X_q \theta=
\sum_{p=1}^6 e_p(\gamma_{pqr}-\gamma_{prq}) X_{\tilde p}
\eeq
are linearly expressed through the operators $X_q$.

Calculating commutators of the operators $X_i$ for the type $[3(21)]$ $h$-spaces
with the help of (\ref{6}),
\bqq{10}
(X_1,X_2)=(X_1,X_4)=(X_2,X_4)=0,
\eeq
$$
(X_1,X_3)=e_2(\gamma_{213}-\gamma_{231})X_2,
\quad
(X_2,X_3)=e_1(\gamma_{123}-\gamma_{132})X_3,
$$
$$
(X_3,X_4)=e_5\gamma_{534}X_4,
\quad
(X_1,X_5)=-e_6\gamma_{651}X_6,
\quad
(X_1,X_6)=-e_5\gamma_{561}X_4,
$$
$$
(X_2,X_5)=-e_6 \gamma_{652} X_6,
\quad
(X_2,X_6)=-e_5 \gamma_{562} X_4,
$$
$$
(X_3,X_5)=e_4 \gamma_{435}X_5-e_6 \gamma_{653} X_6,
\quad
(X_4,X_5)=-e_6 \gamma_{654} X_6,
$$
$$
(X_3,X_6)=e_6 \gamma_{636}X_6-e_5 \gamma_{563} X_4,
\quad
(X_4,X_6)=-e_5 \gamma_{564} X_4,
$$
$$
(X_5,X_6)=-e_5 \gamma_{565}X_4+e_6 \gamma_{656}X_6,
$$
after the change of variables $x^{i'}=\theta^i(x)$, one obtains
\bqq{11}
{\nl{\alpha}\xi}^i=P_{\alpha}(x){\delta_\alpha}^i,
\quad
{\nl{3}\xi}^{\gamma}={\nl{5}\xi}^{\beta}={\nl{6}\xi}^{\beta}={\nl{6}\xi}^5=0,
\eeq
where $\alpha=1, 2, 4$, $\gamma=1, 4, 5, 6$, $\beta=1, 2, 3$.
Here $\theta^i$ are solutions of the complete integrable systems from
(\ref{10}) $X_p \theta=0\;
(p \ne 3)$, $X_q \theta=0\;(q \ne 5)$, $X_r \theta=0\;(r \ne 1)$,
$X_1 \theta=X_2 \theta=X_3 \theta=X_4 \theta=0$,
$X_1 \theta=X_4 \theta=X_5 \theta=X_6 \theta=0$ and
$X_1 \theta=X_2 \theta=X_3 \theta=0$. The first three systems have one
solution each, correspondingly, $\theta^3$, $\theta^5$ and $\theta^1$,
the fourth system has two independent solutions $\theta^5$ and $\theta^6$.
The fifth system has two solutions $\theta^2$ and $\theta^3$, while the sixth system has
three solutions $\theta^4$, $\theta^5$, $\theta^6$.

Using these equalities, from equations (\ref{6}) that do not contain
$\gamma_{pqr}$ one finds
\bqq{12}
\varphi=\frac{3}{2} f_3+c,
\quad
f_i=\lambda_i,
\eeq
$f_1=f_2=f_3(x^3)$ being arbitrary functions of the variable indicated,
$f_4=f_5=f_6=\lambda$ are arbitrary constants.

Integrating the system of equations obtained from (\ref{10}) and taking into account
(\ref{6}), (\ref{11}), (\ref{12}) and also \cite{am2}
$$
g^{ij}=\sum_{h=1}^6 e_h {\nl{h}{\xi}}^i {\nl{\tilde h}{\xi}}^j,
$$
after transformation of coordinates one finds
$$
{\nl{1}{\xi}}^1={\nl{2}{\xi}}^2=1,
\quad
{\nl{3}{\xi}}^2=-\frac{\epsilon x^1}{2A},
\quad
{\nl{3}{\xi}}^3=\frac{1}{2A},
\quad
{\nl{4}{\xi}}^4={\nl{5}{\xi}}^5=(f_3-\lambda)^{-3/2},
$$
$$
g^{44}=e_4(f_3-\lambda)^{-3}(\Sigma+\theta(x^5, x^6)),
\quad
g^{45}=e_4(f_3-\lambda)^{-3},
\quad
g^{66}=e_6(f_3-\lambda)^{-3}.
$$
Calculating components of tensors $g_{ij}$, one obtains
\bqq{13}
g_{ij}dx^idx^j=e_3\lbrace (dx^2)^2+4A dx^1 dx^3+2 \epsilon x^1 dx^2 dx^3+
(\epsilon x^1)^2(dx^3)^2\rbrace+
\eeq
$$
e_4(f_3-\lambda)^3\lbrace2dx^4 dx^5-(\Sigma+\omega)(dx^5)^2\rbrace+
e_6 (f_3-\lambda)^3 (dx^6)^2.
$$
Here
\bqq{14}
A=\epsilon x^2+\theta(x^3),
\quad
\Sigma=3(f_3-\lambda)^{-1},
\eeq
$f_3=\epsilon x^3$, $\epsilon=0, 1$, $\theta(x^3), \omega(x^5, x^6)$ are
functions of their variables, $c, \lambda$ are constants. From the results
obtained, using formulas \cite{am2}
$$
{\nl{h}{\xi}\hspace{-2.5mm}\phantom{a}_{i}}=g_{ij} {\nl{h}\xi}^j,
\quad
a_{ij}=\sum_{h, l=1}^6 {e_h e_l} {\ol a}_{hl}
{\nl{\tilde h}{\xi}\hspace{-2.5mm}\phantom{a}_{i}} {\nl{\tilde
l}{\xi}\hspace{-2.5mm}\phantom{a}_{j}},
$$
one finds
\bqq{15}
a_{ij} dx^i dx^j=f_3 g_{i_1j_1}dx^{i_1}dx^{j_1}+
\eeq
$$
2g_{13}dx^2dx^3+4A \epsilon x^1 (dx^3)^2+
\lambda g_{i_2j_2}dx^{i_2}dx^{j_2}+
g_{45} (dx^5)^2,
$$
\bqq{16}
h_{ij}=a_{ij}+(3\epsilon x^3+c)g_{ij},
\eeq
where $i_1, j_1=1, 2, 3$, $i_2, j_2=4, 5, 6$.

\section{$h$-space of the type $[(32)1]$}

In this case, the commutators are
\bqq{17}
(X_1,X_2)=-e_5 \gamma_{521} X_4,
\quad
(X_1,X_3)=-e_4\gamma_{431}X_5-e_5\gamma_{531}X_4,
\eeq
$$
(X_1,X_4)=-e_3\gamma_{341}X_1,
\quad
(X_1,X_5)=-e_2\gamma_{251}X_2-e_3\gamma_{351}X_1,
$$
$$
(X_1,X_6)=-e_3\gamma_{361}X_1,
\quad
(X_2,X_3)=e_5 \gamma_{523} X_4-e_4 \gamma_{432} X_5-
e_5 \gamma_{532}X_4,
$$
$$
(X_2,X_4)=e_5 \gamma_{524} X_4-e_3 \gamma_{342}X_1,
$$
$$
(X_2,X_5)=e_5 \gamma_{525} X_4-e_2 \gamma_{252}X_2-e_3 \gamma_{352}X_1,
$$
$$
(X_2,X_6)=-e_3\gamma_{362}X_1-e_2\gamma_{262}X_2+e_5\gamma_{526}X_4,
$$
$$
(X_3,X_4)=-e_3\gamma_{343}X_1+e_4\gamma_{434}X_5+e_5\gamma_{534}X_4,
$$
$$
(X_3,X_5)=-e_2 \gamma_{253} X_2+e_4 \gamma_{435} X_5-
e_3 \gamma_{353}X_1+e_5 \gamma_{535} X_4,
$$
$$
(X_3,X_6)=-e_1 \gamma_{163} X_3-e_2 \gamma_{263}X_2-e_3 \gamma_{363}X_1+
e_4 \gamma_{436} X_5+e_5 \gamma_{536}X_4,
$$
$$
(X_4,X_5)=-e_2\gamma_{254}X_2+e_3(\gamma_{345}-\gamma_{354})X_1,
$$
$$
(X_4,X_6)=e_3\gamma_{346}X_1-e_5\gamma_{564}X_4,
$$
$$
(X_5,X_6)=e_2 \gamma_{256} X_2+
e_3 \gamma_{356}X_1-e_4 \gamma_{465} X_5-e_5 \gamma_{565} X_4.
$$
It is clear that the systems $X_p \theta=0\;
(p \ne 3)$, $X_q \theta=0\;(q \ne 6)$ are completely integrable and have one solution
each, $\theta^3$ and $\theta^6$.
The systems $X_1 \theta=X_2 \theta=X_4 \theta=X_6 \theta=0$,
$X_1 \theta=X_4 \theta=X_6 \theta=0$,
$X_1 \theta=X_6 \theta=0$ and the equation $X_6 \theta=0$
are also completely integrable. The first system has solutions $\theta^3$ and $\theta^5$,
the second one has $\theta^2$, $\theta^3$ and $\theta^6$, and
the third one has $\theta^2$, $\theta^3$, $\theta^4$ and
$\theta^5$. The equation has the solution  $\theta^q$ $(q\ne 6)$.
Making the coordinate transformation $x^{i'}=\theta^i(x)$, one obtains in
the new coordinates (with primes omitted)
\bqq{18}
{\nl{s}\xi}^i=P(x){\delta_s}^i,
\quad
{\nl{2}\xi}^{r}={\nl{2}\xi}^5={\nl{3}\xi}^6={\nl{4}\xi}^r=
{\nl{4}\xi}^2={\nl{4}\xi}^5={\nl{5}\xi}^{r}=0,
\eeq
where $s=1, 6$, $r=3, 6$.

From equations (\ref{7}) that do not contain $\gamma_{pqr}$,
one obtains
\bqq{19}
\varphi=\frac{1}{2} f_6+c,
\quad
f_i=\lambda_i,
\eeq
$f_6$ being an arbitrary function of $x^6$,
$f_1=f_2=f_3=f_4=f_5=\lambda$ are arbitrary constants.

Comparing coefficients of same derivatives $\frac{\pl}{\pl x^i}$
in relations (\ref{17}), one obtains the system of $60$ equations for
components of the frame vectors which, after a proper change of variables
leads to
$$
{\nl{1}{\xi}}^1={\nl{2}{\xi}}^2={\nl{3}{\xi}}^3=(f_6-\lambda)^{-1/2},
\quad
{\nl{6}{\xi}}^6=1,
$$
$$
g^{11}=e_3(f_6-\lambda)^{-1}\lbrace(f_6-\lambda)^{-2}+\theta(x^3, x^5)\rbrace,
\quad
g^{12}=e_3(f_6-\lambda)^{-2},
$$
$$
g^{44}=e_5(f_6-\lambda)^{-1}\lbrace(f_6-\lambda)^{-1}+\omega(x^3, x^5)\rbrace,
\quad
g^{45}=e_5(f_6-\lambda)^{-1}.
$$
It follows from above that the tensors $g_{ij}$, $a_{ij}$ and $h_{ij}$ are of
the form
\bqq{20}
g_{ij}dx^idx^j=e_3(f_6-\lambda)
\lbrace (dx^2)^2+2 dx^1 dx^2+\omega(dx^3)^2\rbrace-
\eeq
$$
2e_3 dx^2dx^3+e_4(f_6-\lambda)
\lbrace2dx^4 dx^5-\omega(dx^5)^2\rbrace-e_5(dx^5)^2+
e_6(dx^6)^2,
$$
\bqq{21}
a_{ij} dx^i dx^j=\lambda g_{i_1j_1}dx^{i_1}dx^{j_1}+
2g_{22}dx^2dx^3+
\eeq
$$
g_{23} (dx^3)^2+g_{45} (dx^5)^2+e_6 f_6 (dx^6)^2,
$$
\bqq{22}
h_{ij}=a_{ij}+(f_6+c)g_{ij},
\eeq
where $f_6(x^6), \theta(x^3, x^5), \omega(x^3, x^5)$ are arbitrary functions
of indicated variables, $c, \lambda$ are constants,
$i_1, j_1=1, 2, 3, 4, 5$.

\section{$h$-space of the type $[(321)]$}

In the case of the type $[(321)]$ $h$-space, $\varphi={\rm const}$, hence,
the tensor $h_{ij}$ is covariantly constant. Omitting further calculations,
we write just the final result,
\bqq{23}
g_{ij}dx^idx^j=e_3\lbrace (dx^2)^2+2 dx^1 dx^3-
2 dx^2dx^3+\theta(dx^3)^2\rbrace+
\eeq
$$
e_4\lbrace2dx^4 dx^5-\omega(dx^5)^2\rbrace+e_6(dx^6)^2,
$$
\bqq{24}
a_{ij} dx^i dx^j=\lambda g_{ij}dx^idx^j+
2e_3dx^2dx^3+e_3 (dx^3)^2+e_4 (dx^5)^2,
\eeq
\bqq{25}
h_{ij}=a_{ij}+cg_{ij},
\eeq
where $\theta, \omega$ are arbitrary functions of the variables $x^3, x^5, x^6$,
$c, \lambda$ are arbitrary constants.

\bigskip

The results obtained can be formulated as the following

\bigskip

\noindent
{\bf Theorem 1.} {\it If the tensor $h_{ij}$ of the types $[3(21)]$, $[(32)1]$, $[(321)]$
and the function $\varphi$ satisfy in $V^6$ the Eisenhart equations {\rm (\ref{1})},
then there exists a holonomic coordinate system so that the function $\varphi$ and
the tensors $g_{ij}$, $h_{ij}$ are defined by formulas  {\rm (\ref{12})--(\ref{16})},
{\rm (\ref{19})--(\ref{22})}, {\rm (\ref{23})--(\ref{25})}}.

\section{The quadratic first integrals of the geodesic equations
of the $h$-spaces of the  $[3(21)]$, $[(32)1]$, $[(321)]$ types}

For every solution $h_{ij}$ of equation (\ref{1}) there is a
quadratic first  integral of the geodesic equations (see \cite{am2}),
\bqq{26}
(h_{ij}-4\varphi g_{ij})\dot{x}^i \dot{x}^j={\rm const},
\eeq
where $\dot{x}^i$ is the tangent vector to the geodesic.

Therefore, the quadratic first integrals of the geodesic equations
of the $h$-spaces of the $[3(21)]$, $[(32)1]$, $[(321)]$ types are
determined by formula (\ref{26}), where the tensors $h_{ij}$, $g_{ij}$ and
the function $\varphi$ are obtained in Theorems 1.

\bigskip
I am grateful to professor A.V.Aminova for constant attention to this work
and for useful discussions. The research was partially supported by the RFBR
grant 01-02-17682-a.

\end{document}